\documentclass[11pt,oneside]{amsart}

\usepackage{amscd}
\usepackage{amssymb}
\usepackage[latin1]{inputenc}
\usepackage[pdftex]{graphicx}
\usepackage{verbatim}

\def\N{\mathbb{N}}
\def\Z{\mathbb{Z}}
\def\1{^{-1}}
\def\iff{\Longleftrightarrow}

\def\P{\pi_{1}}
\def\a{\alpha}

\def\f{\phi}

\def\w{\omega}
\def\G{\Gamma}

\def\cal{\mathcal}

\def\el{\'el\'ement}

\newtheorem{lem}{Lemme}[section]
\newtheorem{prop}{Proposition}[section]
\newtheorem{defn}{D\'efinition}
\newtheorem{propdef}{Proposition-D\'efinition}

\title{Groupes à classes de conjugaison  infinies : quelques exemples}

\begin{document}
\maketitle
\begin{center}
{\sc Jean-Philippe PR\' EAUX}\footnote[1]{Centre de recherche de
l'Ecole de l'air, Ecole de l'air, F-13661 Salon de Provence air}\
\footnote[2]{Centre de Math\'ematiques et d'informatique,
Universit\'e de Provence, 39 rue F.Joliot-Curie, F-13453 marseille
cedex 13\\
\indent {\it E-mail :} \ preaux@cmi.univ-mrs.fr}
\end{center}

\begin{abstract}
We consider the group property of being with infinite conjugacy
classes (or {\it icc}, {\it i.e.} $\not= 1$ and of which all
conjugacy classes except 1 are  infinite) in some particular
cases, and more particulary in case of specific
elementary algebraic constructions or of some 3-manifold fundamental groups.\smallskip\\
\noindent\textsc{Résumé.} Nous étudions la propriété pour un
groupe d'être à classes de conjugaison infinis (ou \textit{cci},
{\it i.e.} $\not=1$ et dont toutes les classes de conjugaison
autres que 1 sont infinies) dans certains cas particuliers, et
plus particulièrement dans le cas de constructions algébriques
élémentaires, ou de groupes fondamentaux de certaines 3-variétés.
\end{abstract}

\section{Introduction}

La propriété de groupe d'être à classes de conjugaison infinies
 est particulièrement simple à
énoncer (voir la définition section 2). Elle est motivée par la
caractérisation de Murray et Von-Neumann des algèbres de
Von-Neumann qui sont des facteurs de type II-1 (cf. \cite{MuNe-})
:\smallskip\\
\noindent{\bf Caractérisation de Murray et Von-Neumann.} {\it
Soient $\Gamma$ un groupe et $W_{\lambda}^*(\Gamma)$ son algèbre
de Von-Neumann. Alors $W_{\lambda}^*(\Gamma)$ est un facteur de
type II-1 si et seulement si $\Gamma$ est cci.}

Nous étudions plus particulièrement la stabilité de cette
propriété par certaines constructions algébriques élémentaires, et
le cas des groupes fondamentaux de certaines 3-variétés. L'étude
que nous menons est somme toute élémentaire ;  le lecteur
intéressé pourra trouver des résultats analogues de plus grande
portée dans \cite{cor}, \cite{aogf3v} et \cite{extcci}.\smallskip\\
\indent Ce travail a été effectué durant l'été 2004 alors que
l'auteur se trouvait à l'université de Genève et a été
partiellement financé par le Fonds National Suisse de la Recherche
Scientifique. L'auteur tient à remercier les mathématiciens
Genévois pour leur accueil, et plus particulièrement Goulnara
Arjantseva et Pierre de la Harpe. L'auteur remercie encore Pierre
de la Harpe de lui avoir introduit le problème, ainsi que des
nombreuses discussions et des commentaires apportés.

\section{\bf D\'efinitions et exemples}

\begin{defn}
Un groupe  est dit à classes de conjugaisons infinies (ou cci) si
il est non trivial, et si tout \'el\'ement non trivial a une
classe de conjugaison infinie.
\end{defn}

Il est commode de caract\'eriser cette propri\'et\'e en terme de
centralisateur plut\^ot qu'en terme de classe de conjugaison.

\begin{propdef}
Un groupe non trivial est cci si et seulement si tout \'el\'ement
non trivial a un centralisateur
 d'indice infini.
 \end{propdef}

 \noindent{\bf Preuve.}
 On note $G$ le groupe, et $Z(g)$ le centralisateur d'un \el\ $g\in G$. Soient $g\in G-\{1\}$, et $a,b\in G$.
 Supposons que $aga\1=bgb\1$ :
$$aga\1=bgb\1\iff (b\1 a)g=g(b\1 a) \iff b\1a\in Z(g)$$
Ainsi $aga^{-1}\not= bgb^{-1}$ si et seulement si $a$ et $b$ sont
dans des right-cosets diff\'erents de $G/Z(g)$ et donc, la classe
de conjugaison de $g$, $[[g]]$, est en correspondance bi-univoque
avec $G/Z(g)$. Ainsi $\#[[g]]=\#(G/Z(g))$, ce qui permet de
conclure.\hfill$\square$\\

On peut d\'ej\`a \'etablir quelques exemples \'evidents :

\begin{prop}
\label{obvious}
\begin{itemize}\item[]
\item[(1)] Les groupes libres de rang $\geq 2$ sont cci.
\item[(2)] Les groupes Gromov-hyperboliques non \'el\'ementaires
et sans torsion sont cci.
 \item[(3)] Les groupes ayant un nombre
fini (non nul) d'\el s de torsion (en particulier les groupes
finis), ne sont pas cci. \item[(4)] Les groupes virtuellement
ab\'eliens ne sont pas cci. \item[(5)] Les groupes ayant un centre
non trivial ne sont pas cci.
\end{itemize}
\end{prop}

\noindent{\bf Preuve.} Dans un groupe libre  le centralisateur
d'un \el\ non trivial est cyclique infini. Puisqu'un groupe libre
de rang $\geq 2$ n'est pas virtuellement $\Z$, cela prouve (1).
Dans un groupe hyperbolique le centralisateur d'un \'el\'ement
d'ordre infini est virtuellement cyclique (\cite{cdp}). Ainsi \`a
l'exception des cas o\`u le groupe est fini ou virtuellement
cyclique, le centralisateur d'un \el\ est toujours d'ordre infini,
ce qui prouve (2).
 Pour prouver (3) consid\'erons un
\el\ $g$ de torsion (non trivial). Tous ses conjugu\'es ont m\^eme
ordre que lui, et sont donc aussi de torsion. Ainsi $g$ n'a qu'un
nombre fini de conjugu\'es.
 Pour prouver (4) appelons
$A$ un sous-groupe ab\'elien d'indice fini dans un groupe $G$, et
consid\'erons un \el\ non trivial $g\in A$. Alors $Z(g)\supset A$,
et donc $Z(g)$ est d'indice fini dans $G$.
Pour prouver (5) il suffit de remarquer que tout \el\ du centre a
pour centralisateur le groupe ambient. \hfill$\square$

\section{\bf Stabilité par construction alg\'ebrique}

Afin de consid\'erer le cas d'un produit libre, il nous est
n\'ecessaire d'\'etablir le lemme suivant :

\begin{lem}
Soient $A,B$ des groupes non triviaux. Soit $A$ et $B$ sont tous
deux d'ordre 2, soit $A*B$ contient le groupe libre de rang 2.
\end{lem}

\noindent{\bf Preuve.} Si $A$ et $B$ contiennent chacun un \el\
d'ordre infini $a\in A$, $b \in B$, alors $<a>*<b>$ est un
sous-groupe de $A*B$ isomorphe \`a $F(2)$.

Si $A$ ou $B$ contient un \el\ d'ordre $> 2$, $A*B$ contient un
sous-groupe $<x,y\ |\ x^r,y^s>$, avec $r>2$ et $s\geq 0$,
$s\not=1$. On pose $a=x(xy)x^{-1}$ et $b=yx(xy)x^{-1}y^{-1}$.
Alors avec le th\'eor\`eme d'\'ecriture sous forme r\'eduite, un
mot r\'eduit non trivial sur $a$ et $b$ est non trivial dans
$A*B$, et donc $a$ et $b$ engendrent un groupe libre de rang 2.

Pour finir, si $A$ et $B$ ne sont pas tous deux d'ordre 2, mais ne
contiennent que des \el s d'ordre 2, alors $A*B$ contient $<x,y,z\
|\ x^2,y^2,(xy)^2,z^2>$. Alors comme pr\'ec\'edemment, $a=x(yz)x$
et $b=zx(yz)xz$ engendrent un groupe libre de rang
2.\hfill$\square$

\begin{prop}[\bf Produit libre]
\label{free}
 Soient $A,B$ deux groupes non triviaux. Si $A$ ou $B$
est d'ordre $> 2$, alors $A\ast B$ est cci.
\end{prop}

\noindent {\bf Preuve.} Soit $u$ un \el\ dans un des deux facteurs
; sans perte de g\'en\'eralit\'e on suppose $u\in A$. Avec le
th\'eor\`eme de commutativit\'e dans un produit libre (cf. [MKS]),
le centralisateur de $u$ dans $A*B$ n'est rien d'autre que le
centralisateur de $u$ dans $A$. Puisque $A$ est d'indice infini
dans $A\ast B$, il en va de m\^eme de $Z(u)$. Maintenant si $v$
est un \el\ qui n'est pas dans un des facteurs, son centralisateur
est cyclique infini. Pour conclure il suffit d'utiliser le lemme
pr\'ec\'edent : puisque $A*B$ contient le groupe libre de rang 2,
il ne peut pas \^etre virtuellement $Z$, et donc $A*B$ est cci.\hfill$\square$\\

\noindent{\bf Remarque :} Dans le cas restant, celui du groupe
$\Z_2\ast \Z_2$, ce dernier contient un groupe cyclique infini
d'indice 2 (si l'on note $A*B=<a,b\ |\ a^2,b^2>$, consid\'erer
$<ab>$), et n'est donc pas cci.

\begin{prop}[\bf Sous-groupe d'indice fini]
Si $G$ est cci, et si $K$ est un sous-groupe d'indice fini de $G$,
alors $K$ est aussi cci.
\end{prop}
\noindent{\bf Preuve.} Soit $u\not=1$ dans $K$. On note
respectivement $Z_G(u)$ et et $Z_K(u)$ son centralisateur dans $G$
et $K$ ; trivialement $Z_K(u)=Z_G(u)\cap K$. Puisque $G$ a la
propri\'et\'e H, $Z_G(u)$ est d'indice infini dans $G$. Si
$Z_K(u)$ \'etait d'indice fini dans $K$, il serait aussi d'indice
fini dans $G$, ce qui est impossible puisque $Z_G(u)\supset
Z_K(u)$.\hfill$\square$\\

\noindent{\bf Remarque :} La r\'eciproque est fausse ; il suffit
pour voir cela de consid\'erer le groupe $F_2\times \Z_2$. Elle
devient vraie lorsque l'on suppose que le groupe est sans torsion.

\begin{prop}[\bf Extension finie]
\label{extfini} Si $G$ est sans torsion, et si $K$ est un
sous-groupe d'indice fini de $G$ cci, alors $G$ est aussi cci.
\end{prop}

\noindent{\bf Preuve.} Supposons que $G$ ne soit pas cci, bien
qu'il contienne un sous-groupe d'indice fini cci. Ainsi, il existe
$u\not=1$ dans $G$, tel que $Z_G(u)$ soit d'indice fini dans $G$,
et donc $Z_G(u)\cap K=K_0$ est d'indice fini dans $K$.
N\'ecessairement $u\not\in K$, car cela contredirait le fait que
$K$ est cci.  Mais il est facile de voir qu'il existe $n>1$ tel
que $u^n\in K$. L'ensemble $K_0$ est inclus dans $Z_K(u^n)$, qui
est donc n\'ecessairement d'indice fini dans $K$. Puisque $K$ est
cci, on a $u^n=1$, et donc $G$ a de la
torsion.\hfill$\square$\\

\begin{prop}[\bf Epi-r\'esiduellement cci]
\label{epiresiduel}
 Soit $G$ un groupe ; si pour tout $g\in
G-\{1\}$ il existe un morphisme surjectif $\f$ sur un groupe cci,
tel que $\f(g)\not= 1$ (on dira que $G$ est \'epi-r\'esiduellement
$H$), alors $G$ est cci.
\end{prop}

\noindent {\bf Preuve.} Soit $\f : G \longrightarrow K$ un
\'epimorphisme sur un groupe $K$ cci. Remarquons que puisque $\f$
est surjectif, un conjugu\'e de $\f(g)$ dans $K$ est l'image par
$\f$ d'un conjugu\'e de $g$ dans $G$. Ainsi, puisque $K$ est cci,
si $\f(g)\not=1$ sa classe de conjugaison dans $K$ est infinie, et
donc la classe de conjugaison de $g$ dans $G$ est elle-m\^eme
infinie.\hfill$\square$\\

On en d\'eduit directement :
\begin{prop}[\bf Produit direct]
Soient $A$ et $B$ deux groupes non triviaux. Le produit direct
$A\times B$ est cci si et seulement si $A$ et $B$ sont cci.
\end{prop}

\noindent {\bf Preuve.} Si $A$ ou $B$ est non cci, trivialement
$A\times B$ n'est pas cci ; aussi supposons dans la suite que $A$
et $B$ sont tous deux cci et montrons que $A\times B$ est cci.
Soit $u=ab$, avec $a\in A$ et $b\in B$. Si $u\not= 1$, alors sans
perte de g\'en\'eralit\'e on peut supposer $a\not =1$. En
consid\'erant la projection canonique $A\times B\longrightarrow
A$, la proposition \ref{epiresiduel} permet de
conclure.\hfill$\square$\\

On peut vérifie aussi la stabilité de la propriété d'être cci par
extension.

\begin{prop}[\bf Extension]
Soit $G$ l'extension de deux groupes cci $K$ et $Q$ :
$$1\longrightarrow K\longrightarrow G\overset{p}{\longrightarrow}Q\longrightarrow 1$$
Alors $G$ est aussi cci.
\end{prop}

\noindent{\bf Preuve.} Soit $\w\not=1$ un élément de $G$. Si
$\w\in K$, sa classe de conjugaison dans $K$ est infinie, et il en
est donc de même de sa classe de conjugaison dans $G$. Si
$\w\not\in K$, alors $p(\w)$ a une classe de conjugaison infinie
dans $Q$, et donc puisque $p$ est surjectif, la classe de
conjugaison de $\w$ dans $G$ est aussi infinie.\hfill$\square$

\begin{prop}[\bf Amalgame]
\label{amalgam}
 Si $A$ et $B$ sont deux groupes cci, et si $C=A\cap B$ n'est pas d'indice 1 ou 2 dans
chaque facteur, alors le produit amalgam\'e $A\ *_C B$ est cci.
\end{prop}

\noindent{\bf Preuve.} On peut supposer que $A,B\not=\{1\}$, car
autrement le groupe obtenu n'est autre que $A$ ou $B$. De la
m\^eme fa\c con on peut supposer que $C$ est d'indice $>1$ dans
$A$ et dans $B$.

Soit $u\in A\ *_C B$ ; si $u$ est dans un des facteurs, par
exemple $A$, sa classe de conjugaison dans $A\ *_C B$ contient sa
classe de conjugaison dans $A$, et est donc infinie.

Supposons sans perte de g\'en\'eralit\'e que $C$ est d'indice $>2$
dans $A$. Consid\'erons un \el\  $u\not= 1$ quelconque dans $A\
*_C B$, que l'on peut supposer cycliquemet r\'eduit. Le cas o\`u
$u$ est dans un des facteurs ayant d\'ej\`a \'et\'e trait\'e, on
peut supposer que $u$ s'\'ecrit sous forme cycliquement r\'eduite
$u=a_1b_1.\cdots .a_nb_n$, avec $n\geq 1$. L'ensemble des
right-cosets $A/C$ contient au moins trois \el s distincts $1.C$,
$\a.C$ et $\a^2.C$. On choisit $a\in A$ hors de la classe $1.C$ et
dans une classe diff\'erente de celle de $a_1\1$. On choisit $b\in
B-C$. Avec ce choix, tous les \el s $(ba)^n.u.(ba)^{-n}$ de la
classe de conjugaison de $u$ sont distincts. Ainsi tout \el\ non
trivial de $A\ *_C B$ a une classe de conjugaison infinie, ce qui
conclut la preuve.\hfill$\square$\\

Le m\^eme mod\`ele de preuve s'applique pour montrer :
\begin{prop}[\bf Extension HNN]
\label{hnn}
 Si $A$ est cci, et si $C$ est
un sous-groupe d'indice $>2$ dans $A$, l'extension HNN $A*_C$ est
cci.
\end{prop}

\noindent{\bf Preuve.} Soit $A*_\f$, avec $\f:C\longrightarrow
\f(C)$. On note $t$ une lettre stable, telle que $tct\1=\f(t)$
pour tout $c\in C$. Comme pr\'ec\'edemment on peut se ramener
facilement au cas d'un \el\ cycliquement r\'eduit de longueur
$>1$, $u=a_1t^{n_1}.\ldots.\a_pt^{n_p}$, avec $n_p\not=0$. On a au
moins 3 right-cosets dans  $A/C$ : $1.C$, $\a.C$, et $\a^2.C$. On
choisit $a\in A$ hors de $1.C$ et de la classe de $a_1\1$. Alors
tous les \el s $(ta)^n u (ta)^{-n}$ pour $n>0$ sont distincts dans
la classe de conjugaison de $u$, ce qui permet de conclure.\hfill$\square$\\

Et on obtient en corollaire imm\'ediat :
\begin{prop}[\bf Graphes de groupes]
Si $\G=\P (\cal{G},X)$ est le groupe d'un graphe de groupe dont
tous les groupes de sommet sont cci, et tel que pour toute ar\^ete
$e$ d'origine $s_0$ et d'extr\'emit\'e $s_1$, $G(e)$ est d'indice
$>2$ dans $G(s_0)$ ou bien dans $G(s_1)$, alors $\G$ est cci.
\end{prop}

\section{\bf Groupes de 3-vari\'et\'es}

Nous considérons d'abord le cas des groupes fondamentaux de fibrés
de Seifert (cf.\cite{seifert}). En utilisant le fait que dans un
groupe de Seifert, la classe d'une fibre r\'eguli\`ere a un
centralisateur d'indice fini (cf. \cite{js}), on \'etablit
imm\'ediatement :
\begin{prop}[\bf Fibr\'es de Seifert]\label{seifert}
Les groupes de fibr\'es de Seifert ne sont pas cci.
\end{prop}
\vskip 0.2cm

Consid\'erons maintenant le cas des vari\'et\'es hyperboliques de
volume fini.  Il montre une autre application directe de la
proposition \ref{epiresiduel}.

\begin{prop}[\bf Hyperboliques et $\chi=0$]
Les groupes de 3-vari\'et\'es hyperboliques de volume fini non
\'el\'emen\-taires sont cci.
\end{prop}

\noindent{\bf Preuve.} Consid\'erons une vari\'et\'e $M$
v\'erifiant les hypoth\`eses : $M$ est une 3-vari\'et\'e
hyperbolique \`a bord torique ou vide, qui n'est ni un tore
solide, ni un tore \'epaissi. En outre, la classification des
isométries hyperboliques montre que $\pi_1(M)$ est sans torsion
(cf. \cite{bp}), et donc avec la proposition \ref{extfini}, en
considérant le revêtement d'orientation de $M$ on supposera $M$
orientable.

Le th\'eor\`eme de chirurgie hyperbolique de Thurston affirme
qu'il existe une suite de vari\'et\'es hyperboliques ferm\'ees
$(M_n)_n$, obtenues par obturation de Dehn sur $M$, convergeant
vers $M$ pour la topologie g\'eom\'etrique (cf. \cite{bp}). Ainsi
on a une suite d'\'epimorphismes $(\f_n: \P M\longrightarrow\P
M_n)_n$. Par d\'efinition de la topologie g\'eom\'etrique, si
$u\not=1\in \P M$, pour $n\gg 0$, $\f_n(u)\not=1$ ; et tous les
groupes $\P M_n$ sont hyperboliques au sens de Gromov. De plus
pour $n\gg 0$ ils sont non \'el\'ementaires (ceci car  leurs
points fixes dans $\partial \mathbb{H}^3$ convergent vers les
points fixes de $\P M$ et sont donc pour $n\gg 0$ en nombre $>2$).
Ainsi avec la proposition \ref{obvious} (2), $\P M$ est
\'epi-r\'esiduellement cci, et la proposition \ref{epiresiduel}
permet de conclure.\hfill
$\square$\\

Nous traitons maintenant le cas des groupes de 3-variétés
orientables suffisament large. Rappelons qu'un 3-variété
orientable est dite suffisament large si elle contient une surface
incompressible à deux faces ; une 3-variété orientable est Haken
si elle est irréductible et suffisament large.

\begin{prop}[{\bf Le cas Haken orientable}]\label{haken}
Soit $M$ un 3-variété Haken orientable ; alors soit $\pi_1(M)$ est
cci soit $M$ est un fibré de Seifert.
\end{prop}

\noindent{\bf Preuve.} Nous avons vu qu'un fibré de Seifert a un
groupe non cci ; aussi supposons que $\pi_1(M)$ est non cci et
montrons que $M$ est un fibré de Seifert. Rappelons que le groupe
d'une 3-variété Haken est non trivial (car il contient un groupe
de surface infini) ; ainsi, soit $u\not=1$ dans $\pi_1(M)$ ayant
un centralisateur $Z(u)$ d'indice fini, et soit $p:
N\longrightarrow M$ le revêtement fini associé à $Z(u)$.
Clairement $\pi_1(N)=Z(u)$ a un centre non trivial. De plus
nécessairement $N$ est Haken : le corollaire 13.5 de \cite{hempel}
montre que $N$ est irréductible, et de plus si $S$ est une surface
incompressible à deux faces dans $M$, $p^{-1}(S)$ est une surface
incompressible à deux faces dans $N$. Le corollaire II.5.4 de
\cite{js} montre que $N$ est nécessairement un fibré de Seifert et
le théorème II.5.3 (\cite{js}) montre qu'il en est alors de même
pour $M$.\hfill $\square$

\begin{prop}[{\bf Le cas orientable suffisament large}]
Soit $M$ une 3-variété orientable et suffisament large ayant un
groupe fondamental non diédral infini ; alors soit $\pi_1(M)$ est
cci, soit la complétée de Poincaré $\cal P(M)$ est un fibré de
Seifert.
\end{prop}

\noindent {\bf Preuve.} Considérons la décomposition de
Kneser-Milnor (cf. \cite{hempel}) de $M$ :
$$M=M_1\#\cdots\# M_n\#B_1\#\cdots \# B_p\# C_1\#\cdots C_q$$
avec $n,p,q\in\N$, et où les $M_i$ sont non simplement connexes et
premiers, les $B_i$ sont des boules et les $C_i$ des fausses
3-sphères. Puisque $M$ est suffisament large, $\pi_1(M)$ est non
trivial (car il contient un groupe de surface infini) et donc
$n>0$. La complétée de Poincaré $\cal P(M)$ est alors définie
comme étant :
$$\cal P(M)=M_1\#\cdots\# M_n$$
et de plus $\pi_1(M)=\pi_1(\cal P(M))=\pi_1(M_1)*\cdots
*\pi_1(M_n)$. Ainsi avec la proposition \ref{seifert} si $\cal
P(M)$ est un fibré de Seifert, alors $\pi_1(M)$ n'est pas cci.
Supposons dans la suite que $\pi_1(M)$ n'est pas cci. Alors avec
la proposition \ref{free} soit $\pi_1(M)$ est diédral infini, soit
$n=1$. Ainsi sous nos hypothèses $\cal P(M)=M_1$ et soit
$M_1\equiv S^2\times S^1$ soit $M_1$ est irréductible. Dans le
premier cas $\cal P(M)$ est un fibré de Seifert, et dans le second
cas, puisque $M$ contient une surface incompressible à deux faces,
$M_1$ est Haken, et la proposition \ref{haken} donne la
conclusion.\hfill$\square$\\

Dans \cite{aogf3v} ces résultats ont été généralisés au cas de
groupes fondamentaux de 3-variétés quelconques.

\vskip 0.4cm


\end{document}